\documentclass[10pt]{amsart}

\usepackage{amsfonts}

\usepackage{graphicx}
\usepackage{amsmath}
\usepackage{setspace}
\usepackage{amsthm, amssymb,enumerate,color, bbm, version}
\usepackage[left=3cm,right=3cm,top=2.5cm,bottom=3cm]{geometry}

\newcommand{\halmos}{{\mbox{\, \vspace{3mm}}} \hfill
\mbox{$\Box$}}

\RequirePackage[colorlinks]{hyperref}

\newtheorem{thm}{Theorem}[section]
\newtheorem{lem}[thm]{Lemma}
\newtheorem{cor}[thm]{Corollary}
\newtheorem{prop}[thm]{Proposition}

\theoremstyle{remark}
\newtheorem{rem}{Remark}

\newcommand{\Exp}{{\mathbb E}}

\newcommand{\eqdistr}{\stackrel{\scriptstyle  d}{=}}

\newcommand{\convweak}{\stackrel{\rm w}{\to}}

\newcommand{\dto}{\xrightarrow{d}}

\newcommand{\fidi}{\xrightarrow{\text{fidi}}}
\newcommand{\dsum}{\displaystyle\sum}

\newcommand{\toi}{{\to \infty}}
\newcommand{\eind}{\stackrel{d}{=}}
\newcommand{\1}{\mathbb I}
\newcommand{\bbr}{{\Bbb R}}

\newcommand{\bbn}{{\Bbb N}}

\newcommand{\law}{\mathcal{L}}

\title{Heavy tailed branching process with immigration}

\author{Bojan Basrak}
\address{Department of Mathematics, University of Zagreb, Bijeni\v{c}ka 30, 10000 Zagreb, Croatia.} \email{bbasrak@math.hr}

\author{Rafa{\l} Kulik}
\address{Department of Mathematics and Statistics, University of Ottawa, K1N 6N5 Ottawa, Canada.}
\email{rkulik@uottawa.ca}

\author{Zbigniew Palmowski}
\address{Mathematical Institute, University of Wroc{\l}aw, Poland.}
\email{zbigniew.palmowski@gmail.com}

\thanks{Research of the first author was partially supported by
the Croatian government under the grant MZOS 037-0372790-2800.
Research of the second author was partially supported by NSERC
Grant. Research of the third author was partially supported by the
Ministry of Science and Higher Education of Poland under the grant
N201 412239.}

\date{\today}
\subjclass[2010]{60G51, 60G50, 60K25} %
\begin{document}
\begin{abstract}
In this paper we analyze a branching process with immigration
defined recursively by $X_t=\theta_t\circ X_{t-1}+B_t$ for a
sequence $(B_t)$ of i.i.d. random variables and random mappings $
\theta_t\circ x:=\theta_t(x)=\sum_{i=1}^xA_i^{(t)} \; , $ with
$(A_i^{(t)})_{i\in \mathbb{N}_0}$ being a sequence of
$\mathbb{N}_0$-valued i.i.d. random variables independent of $B_t$.
We assume that one of generic variables $A$ and $B$ has a regularly
varying tail distribution. We identify the tail behaviour of the
distribution of the stationary solution $X_t$. We also prove CLT for
the partial sums that could be further generalized to FCLT. Finally,
we also show that partial maxima have a Fr\'echet limiting
distribution.
\end{abstract}
\maketitle
 \noindent {\sc Keywords:} Branching process with
immigration $\star$ regularly varying distribution $\star$ INAR

\noindent {\sc Short title:} Heavy tailed branching process

\pagestyle{myheadings} \markboth{\sc B.\ Basrak --- R. Kulik
--- Z.\ Palmowski} {\sc Heavy tailed branching process}

\tableofcontents

%\newpage

\section{Introduction}
Branching processes with immigration are one of the corner stone
models in applied probability and Markov chain theory on a countable
state spaces. They were studied extensively in the literature, the
classical references are: \cite{General1, General2, General3,
General4, General5, General6, General7, General8, General9,
General10}. The process considered in this paper appears also in
various models in queueing theory (see \cite{Queueininne1,
Queueininne2, Queueininne4, Queueininne5}), polling systems (with
possible non-zero switch-over times; see \cite{Queueinpolling1,
Queueinpolling2, Queueinpolling3}), infinite server queues (e.g.
\cite{Queueininne3}) and processor sharing queues (e.g.
\cite{Queueingps}). This process can be also used as a model of
packet forwarding in delay-tolerant mobile ad-hoc networks (see
\cite{Queueininne3} for details). In the time series theory, the
model considered here contains a class of so-called INteger
AutoRegressive (INAR) processes, see \cite{Biology1}. Some possible
applications of these processes, starting with medicine and
biological sciences are discussed in \cite{Biology1b, Biology2}. The
INAR processes were originally introduced as a discrete counterpart
of AR processes. However,
their nonlinear structure represents a challenge for theoretical analysis.\\

We study a process $(X_t)_{t\in \mathbb{N}_0}$ indexed over
$\mathbb{N}_0=\{0,1,2,\ldots\}$. Here, $X_0$ is a random
variable independent of an i.i.d. sequence of random mappings
$(\theta_t)_{t\in \mathbb{Z}}$, $\theta_t:\mathbb{N}_0\to
\mathbb{N}_0$ defined as follows: for each $t\in \mathbb{Z}$ and
for a sequence $(A_i^{(t)})_{i\in \mathbb{N}_0}$ of
$\mathbb{N}_0$-valued i.i.d. random variables, we have:
$$
\theta_t\circ x:=\theta_t(x)=\sum_{i=1}^xA_i^{(t)} \; ,\qquad x\in\mathbb{N}_0.
$$
That is, $\theta_t$ maps an integer $x$ into a random integer with
an interpretation that each of $x$ individuals in the $(t-1)$th
generation leaves behind a random number of children, and all
these numbers are independent and have the same distribution as some
generic random variable, say $A$. To introduce immigration in the
model, we assume that another i.i.d. sequence $(B,B_t\,,{t\in
\mathbb{Z}})$ of $\mathbb{N}_0$-valued random variables is given
independently of the sequence $(\theta_t)_{t\in \mathbb{Z}}$.
Then $(X_t)_{t\in \mathbb{N}_0}$ satisfies:
\begin{equation}\label{eq:defINAR}
X_t=\theta_t\circ X_{t-1}+B_t\qquad \mbox{\rm for each } t\ge 1\; ,
\end{equation}
or in an alternative notation
$$
X_t=\sum_{i=1}^{X_{t-1}}A_i^{(t)}+B_t\qquad \mbox{\rm for each }
t\ge 1\; .
$$
Observe that the random mappings $\theta_t$ by definition
satisfy
$$
 \theta_t \circ (x+y) =  \theta^{(1)}_t \circ x+ \theta^{(2)}_t \circ y\;,
$$
where $\theta^{(i)}_t$, $i=1,2$ on the right-hand side are independent
with the same distribution as $\theta_t$.\\

Denote by $f$ and $g$ the probability generating functions of $A$ and $B$, respectively, i.e.
$$
f(z)=E(z^A)\;, \qquad g(z)=E(z^B).
$$
Following Foster and Williamson \cite{fosterwill}, the Markov chain $(X_t)$ is ergodic with unique stationary distribution if and only if
$$
\int_0^1 \frac{1-g(s)}{f(s)-s}ds <\infty .
$$
In terms of moments, sufficient conditions are given in
\cite{Seneta}. If
\begin{equation}\label{basicass1}
0<\mu:=E(A)<1\qquad \mbox{and}\qquad E(\ln (1+B))<\infty\;,
\end{equation}
then the chain is ergodic with unique stationary distribution and even strongly mixing.
From now on we assume that (\ref{basicass1}) holds.

The unique stationary representation of $(X_t)_{t\in\mathbb{N}_0}$ is given by
\begin{equation}\label{eq:stat-representation}
X_t\eqdistr B_t+\sum_{k=1}^{\infty}\theta^{(t-k)}_t\circ \cdots \circ
\theta^{(t-k)}_{t-k+1} ( B_{t-k}) =: B_t+\sum_{k=1}^{\infty}
\bigotimes_{i=0}^{k-1} \theta^{(t-k)}_{t-i} (B_{t-k})=: \sum_{k=0}^{\infty}
C_{t,k} \; .
\end{equation}
Since $(\theta^{(t)}_i)$ are i.i.d., it follows that
$(C_{t,k})_{k\in\mathbb{N}_0}$ is another sequence of independent
random variables.\\

In \cite{Seneta} one can find corresponding representation of the
stationary distribution in terms o probability generating functions.
For completeness, we show in Lemma~\ref{LemaSt} that the random
series above converges with probability 1. From there,  it is
straightforward to check that such $X_t$ satisfy \eqref{eq:defINAR}.

We will also consider a Markov
chain $(X'_t)_{t\in\mathbb{N}_0}$ which evolves as
$$
X'_t= \max\{ \theta_t\circ X'_{t-1}, B_t\} \qquad \mbox{\rm for each } t\ge 1\; .
$$
The unique stationary distribution of $(X_t')$ exists since $X_t'\leq X_t$ and it is given by:
\begin{equation}\label{eq:max-process}
X'_t\eqdistr \sup\{ B_t,\theta^{(t-k)}_t\circ \cdots \circ
\theta^{(t-k)}_{t-k+1} ( B_{t-k}) : k=1,2,\ldots \}  =\sup
\{ C_{t,k} : k\geq 0 \} \; .
\end{equation}
Since $X'_t$  are dominated by $X_t$ in
\eqref{eq:stat-representation}, the supremum above is  a.s. finite.
Direct calculation shows that $(X'_t)$ define strictly stationary
sequence.

In this paper we identify the tail behaviour of the distribution of
the stationary solution $X_t$ under assumption that the generic size
of immigration $B$ or generic size of offsprings $A$ has regularly
varying distribution. Although there is an extensive literature on
behaviour of heavy tailed compound sums (see e.g. \cite{dima} and
reference therein) or heavy tailed random difference equations (see
\cite{HultSamo} and references therein) those results do not seem to
produce in a straightforward manner the asymptotics of stationary
distribution for branching processes. Surprisingly, literature on
tail asymptotics for branching processes is limited; see e.g.
Corollary 2 in \cite{dima}.

We also prove CLT for the heavy tailed partial sums that could be
further generalized to FCLT (see e.g. \cite{limit1}). Furthermore,
we show that partial maxima have Fr\'echet limiting distribution.
%We also state finite variance CLT for completeness. \\

The paper is organized as follows.
%In Section \ref{Stationary solution} we give representation of stationary
%distributions of $X_t$ and $X_t'$ and identify sufficient conditions for their existence.
In Section \ref{Tail behaviour} we find the tail behaviour of
stationary sequences $(X_t)_{t\in\mathbb{N}_0}$ and
$(X_t')_{t\in\mathbb{N}_0}$ under two different regimes. We use this
result in establishing CLT for the process $(X_t)$ %and FCLT for both processes which is described
in Section \ref{Asymptotics for sums and maxima}.

%\section{Stationary solution}\label{Stationary solution}

\section{Tail behaviour}\label{Tail behaviour}

\subsection{Regularly varying immigration (Model I)}\label{section:rv-immigrations}
We will assume that $A$ and $B$ satisfy the following conditions:
\begin{equation}\label{eq:1}
0< \mu=E(A)<1 \; ,
\end{equation}
\begin{equation}\label{eq:2}
P(B>x)=x^{-\alpha}L(x) \; ,
\end{equation}
for some $\alpha\in (0,2)$ and a slowly varying function $L$.
We consider here the case $\alpha\in (0,2)$ only. For $\alpha>2$ proofs become much more involved, however, a technique is clearly suggested by the case $\alpha\in [1,2)$, see  Remark \ref{rem:alpha3} below.
For $\alpha\in [1,2)$, we also assume that
\begin{equation}\label{eq:3a}
E(A^2)<\infty.
\end{equation}
In particular, it means that the tail of $A$ is lighter than that of $B$ in the sense that
\begin{equation}\label{eq:3}
P(A>x)={\rm }o(P(B>x)) \; .
\end{equation}
The conditions above are needed, in particular, to conclude that
\begin{equation}\label{eq:compound_sum_asymptotics}
P\left(\sum_{i=1}^B A_i>x\right) \sim P(B>x/\mu),
\end{equation}
 where $A_i$ are i.i.d. copies of $A$, independent of $B$.
To establish (\ref{eq:compound_sum_asymptotics}), in \cite{Faetal06}
the authors assume (\ref{eq:1})-(\ref{eq:2}) and (\ref{eq:3})
together with $E(B^{\alpha})<\infty$ if $\alpha=1$; see Prop. 4.3 in
\cite{Faetal06}. On the other hand, in \cite[Theorem 3.1]{RoSe} the
authors assume that $B$ is consistently varying (which is implied by
(\ref{eq:2})), $E(A^r)<\infty$ for some $r>1$. Furthermore, they
assume that either
\begin{equation}\label{eq:BS_1}
E(B)<\infty\;\;\; \text{together with  } P(A>x)={\rm o}(P(B>x))
\end{equation}
or
\begin{equation}\label{eq:BS_2}
E(B)=\infty \;\;\; \text{together with  }  xP(A>x)={\rm o}(P(B>x)).
\end{equation}
We note that the latter follows from (\ref{eq:3a}) and the fact that mean of $B$ is infinite.

In the sequel it is useful to introduce the following
random variables
\begin{equation} \label{e:tildeA}
 \tilde A^{(k)} = \theta_k\circ \cdots \circ \theta_{1} \circ 1
\end{equation}
and assume that for each $k$ an i.i.d. sequence $(\tilde
A^{(k)}_i)_{i\in\mathbb{N}_0}$ is given  with the same distribution
as $\tilde A^{(k)}$. Note also that $\tilde A^{(1)}_i \eind A$. In
general we have
$$
 \tilde A^{(k)}_i \eind \sum_{j=1}^A  {\tilde A^{(k-1)}_j}
$$
with all the random variables on the right hand side being
independent. From the above identity and $E(A)=\mu$ we have
$E(\tilde A_i^{(k)})=\mu^k$. Moreover, for $C_{t,k}$ defined in
(\ref{eq:stat-representation}) it holds
\begin{equation}\label{eq:CtildeA}
 C_{t,k} \eind  \sum_{j=1}^{B_{t-k}}  {\tilde A^{(k)}_j}.
\end{equation}

We start our analysis by showing that under our basic assumption
(\ref{basicass1})
 the random series in \eqref{eq:stat-representation} converges a.s.
 and hence $X_t$ in the same expression are properly defined.
\begin{lem}\label{LemaSt}
 Assume that \eqref{basicass1} holds. Then
 $\sum_{k=0}^{\infty}
C_{t,k} $ converges a.s.
\end{lem}
\noindent {\it Proof.}
Take an arbitrary $\varepsilon>0$ such that $\mu+\varepsilon<1$.
Observe that
\begin{eqnarray*}
\lefteqn{P(C_{t,k} \not = 0) = P \left(\sum_{j=1}^{B_{t-k}}  {\tilde A^{(k)}_j} \not = 0\right) }\\
& \leq & P\left( B \geq (\mu+\varepsilon)^{-k}\right) +
P \left(  \dsum_{j=1}^{\lfloor(\mu+\varepsilon)^{-k}\rfloor}  {\tilde A^{(k)}_j} >\varepsilon \right)\\
& \leq & P (\ln (1+B) \geq - k \ln (\mu + \varepsilon) )
+
\dfrac{\mu^k}{(\mu-\varepsilon)^k \varepsilon},\\
\end{eqnarray*}
where the last expression follow by the Markov inequality and
$E(\tilde A_i^{(k)})=\mu^k$. Since $E \ln (1+B)<\infty$, the sum of
probabilities above for $k=1$ to $\infty$
 is finite. Hence Borel-Cantelli lemma shows that with probability 1 only finitely
 many $C_{t,k}$'s in \eqref{eq:stat-representation} are different from 0.
\halmos

Condition~(\ref{eq:2}) ensures that the upper tail of the immigration distribution
is regularly varying.
We first show that this property is inherited by the random variables $X_t$
and $X'_t$
introduced in \eqref{eq:stat-representation} and \eqref{eq:max-process}.

\begin{thm}\label{eq:Theorem-1}
Under the conditions (\ref{basicass1}), (\ref{eq:1})-(\ref{eq:3a})
and (\ref{eq:BS_1}) or (\ref{eq:BS_2}), we have
\begin{equation}\label{eq:tail-behaviour}
\lim_{x\to\infty} \frac{P(X_t>x)}{P(B>x)}
= \lim_{x\to\infty} \frac{P(X'_t>x)}{P(B>x)}
 =\sum_{k=0}^{\infty}
\mu^{k\alpha}\; .
\end{equation}
\end{thm}
\begin{rem}{\rm
Our model is $X_t=\sum_{i=1}^{X_{t-1}}A_i^{(t)}+B_t$. If we assume
that random variables $A_i^{(t)}$ are equal to $\mu\in (0,1)$, then
we obtain $X_t=\mu X_{t-1}+B_t$, that is the AR(1) model. For such
models it is well-known that
$$P(X_t>x)\sim \sum_{k=1}^{\infty}\mu^{k\alpha}P(B>x). $$
Hence, the tail of the branching model with random offspring
$A_i^{(t)}$ is the same as that of the averaged AR(1) model.
 }
\end{rem}
\subsection{Proof of Theorem \ref{eq:Theorem-1}}

For $r,k =1,2,\ldots$ denote
$$
  m_r(k):=E(\tilde A^{(k)}_1)^r=E\left(\sum_{i=1}^A \tilde A_i^{(k-1)}\right)^r,
$$
with a convention $\tilde A_i^{(0)}=A_i$, where $A_i$ are i.i.d. with the same distribution as $A$.
The following lemma turns out to be very useful in the sequel

\begin{lem} \label{lem:bounds}
For $k\geq 2$ we have:
\begin{equation}\label{eq:first-bound}
 E\left(\tilde A_1^{(k)}\right)^2  \le E(A^2) (k+1)\mu^k .
\end{equation}
Moreover, if $E(B^2)<\infty$, then
$$
E\left(\sum_{j=1}^B \tilde A_j^{(k)}\right)^2\le  E(B)
E(A^2)(k+1)\mu^k+E(B^2)\mu^{2k} .
$$
\end{lem}
\noindent {\it Proof.} The proof of (\ref{eq:first-bound}) is by
induction on $k$. For $k=2$ we bound
$$
E\left(\sum_{j=1}^A A_j\right)^2,
$$
where all random variables $A,A_j$ are i.i.d. Simple conditioning argument yields
the bound
$$
E(A) E(A^2) +(E(A))^2E(A^2) =E(A^2)\mu (1+\mu)\le 2 E(A^2) \mu
$$
that matches (\ref{eq:first-bound}) for $k=2$.

For $k\geq 3$ we have
\begin{align*}
m_2(k)&= \sum_n E\left(\sum_{j=1}^n \tilde A_j^{(k-1)}\right)^2 P(A=n) = \\
&=  \sum_n P(A=n)\sum_{j=1}^n E\left(\tilde A_j^{(k-1)}\right)^2 \\
&+  \sum_n P(A=n) \sum_{j,l\atop j\not=l}
E\left(\tilde A_j^{(k-1)}\right)E\left(\tilde A_l^{(k-1)}\right)\\
&\le E(A) E(A^2) k\mu^{k-1}+E(A^2) \mu^{2k} \le E(A^2) (k+1)\mu^k.
\end{align*}
%Computation without induction is as follows:
%\begin{align*}
%m_2(k) &\le  \mu E(\tilde A^{(k-1)})^2 +E(A^2) \mu^{2k} \\
%&= \mu^2 E(\tilde A^{(k-2)})^2 +\mu E(A^2)\mu^{2(k-1)}+E(A^2)\mu^{2k} \\
%&= \mu^3 E(\tilde A^{(k-3)})^2 +\mu^2 \mu^{2(k-2)}E(A^2)+\mu \mu^{2(k-1)}E(A^2)+E(A^2)\mu^{2k}\\
%&= \cdots \\
%&=\mu^k E(\tilde A^{(0)})^2 +E(A^2)\mu^k\left\{\mu+\mu^2+\cdots+\mu^k\right\}\le E(A^2)(k+1)\mu^k .
%\end{align*}
Furthermore,
\begin{align*}
E\left(\sum_{j=1}^B \tilde A_j^{(k)}\right)^2 & \le
\sum_{n=1}^{\infty} nP(B=n) E(\tilde A_0^{(k)})^2 +
\sum_{n=1}^{\infty}n^2 P(B=n)E(\tilde A_0^{(k)})^2\\
&\le E(B) E(A^2)(k+1)\mu^k+E(B^2)\mu^{2k}.
\end{align*}
This completes the proof of Lemma \ref{lem:bounds}. \hfill
$\Box$

\noindent {\it Proof of Theorem \ref{eq:Theorem-1}.} \noindent

\noindent {\bf Step 1:} Large deviation results, such as those
developed in \cite{Na79}, \cite{Faetal06} (see also \cite{RoSe})
yield that under the conditions (\ref{eq:1})-(\ref{eq:3a}) and
either (\ref{eq:BS_1}) or (\ref{eq:BS_2}) we have,
$$ P(\theta_t(B_{t-1})>x) \sim
P(B>x/\mu)\;.
$$
By induction,
\begin{equation}\label{eq:5}
P(C_{t,k}>x)=P\left(\bigotimes_{i=0}^{k-1} \theta_{t-i}^{(t-k)} (B_{t-k})>x
\right)\sim P(B>x/\mu^k)\; .
\end{equation}
Thus, each of $C_{t,k},\, k\geq 0$, is regularly varying and since
they are independent, the random variables
$$
X_{t,m} = \sum_{k=0}^{m} C_{t,k}
$$
and
$$
X'_{t,m} = \max\{ C_{t,k} : {k=0,\ldots,m}\}\,
$$
both have the same property and satisfy
\begin{equation}\label{eq:4}
\lim_{x\to\infty}
\frac{P(X_{t,m}>x)}{P(B>x)}
=
\lim_{x\to\infty}
\frac{P(X'_{t,m}>x)}{P(B>x)}
=\sum_{k=0}^{m}\mu^{k\alpha}\,;
\end{equation}
see \cite{Re07} for instance.

\noindent {\bf Step 2:} To bound $\liminf$ in (\ref{eq:tail-behaviour}) from below
note that:
\begin{eqnarray*}
\liminf_{x\to\infty}\frac{P(\sum_{k=0}^{\infty}C_{t,k}>x)}{P(B>x)} \ge
\lim_{x\to\infty}\frac{P(\sum_{k=0}^{m}C_{t,k}>x)}{P(B>x)} \; .
\end{eqnarray*}
Therefore,
$$
\lim_{x\to\infty}\frac{P(\sum_{k=0}^{\infty}C_{t,k}>x)}{P(B>x)} \ge
\sum_{k=0}^{\infty}\mu^{k\alpha}\; ,
$$
by applying (\ref{eq:4}) and letting $m\to\infty$.

\noindent {\bf Step 3:} To establish an upper bound for the tail of $X_t$ in (\ref{eq:tail-behaviour}) it is enough to show that
\begin{equation}\label{eq:upper_bound_Xt}
\lim_{k_0\to\infty} \limsup_{x\to\infty}\frac{P(\sum_{k =
k_0}^{\infty}C_{t,k}>x)}{P(B>x)}=0.
\end{equation}
Similarly, to obtain an upper bound for the tail of $X'_t$, it is enough to show
that
\begin{equation}\label{eq:upper_bound_Xtprime}
\lim_{k_0\to\infty} \limsup_{x\to\infty}\frac{P(\sup_{k \geq
k_0}C_{t,k}>x)}{P(B>x)}=0.
\end{equation}
However, once we show \eqref{eq:upper_bound_Xt}, (\ref{eq:upper_bound_Xtprime}) follows immediately.

Observe that:
\begin{eqnarray}
\lefteqn{\frac{P(\sum_{k = k_0}^{\infty} C_{t,k}>x)}{P(B>x)}}\label{eq:tail4sumofCs} \\
& \leq & \frac{P(\sup_{k \geq k_0} B_{t-k}>x(1-\varepsilon)/\mu^k)}{P(B>x)}
+ \frac{P(\sum_{k = k_0}^{\infty} C_{t,k} \1_{\{B_{t-k}< x(1-\varepsilon)/\mu^k \}} >x)}{P(B>x)}. \nonumber
\end{eqnarray}
The first term on the right hand side is bounded above by:
$$
 \sum_{k = k_0}^{\infty} \frac{P(B>x(1-\varepsilon)/\mu^k)}{P(B>x)}\,.
$$
One can use the Potter's bounds (see \cite{Re07}) now
 to see that its limit is zero if we let first $x$ and then $k_0$
converge to $\infty$. The second term is more difficult to handle,
hence we split the proof in different cases with respect to the
value of $\alpha$.

\noindent{\textbf{Case}} $0<\alpha<1$: For the second term in
\eqref{eq:tail4sumofCs} we apply the Markov inequality to bound it
by
$$
\sum_{k = k_0}^{\infty}  \mu^k \frac{  E\left( B_{t-k}
\1_{\{B_{t-k}< x(1-\varepsilon)/\mu^k \}}\right)}{x P(B>x)}.
$$
By Karamata's theorem in combination with the Potter's bounds each summand above is
bounded by
$$
 (1+\varepsilon) \mu^k  \frac{1}{1-\alpha}
 \left(\frac{1-\varepsilon}{\mu^k}\right)^{1-\alpha}\,,
$$
for $x,\, k_0$ large enough.
Hence, we observe that (\ref{eq:upper_bound_Xt}) holds once we take $\lim_{k_0\to\infty} \limsup_{x\to\infty}$ in \eqref{eq:tail4sumofCs}.

Note that a nonnegative random variable $Y$ is regularly varying with index $\alpha$ if and only if
$Y^j,\ j>0$ is regularly varying with index $\alpha/j$. This remark turns out to be very useful
in the rest of the proof.

\noindent{\textbf{Case}} $1\leq \alpha<2$:
Note that (\ref{eq:upper_bound_Xt}) is equivalent to a requirement:
$$\lim_{k_0\to\infty} \limsup_{x\to\infty}\frac{P(\sum_{k =
k_0}^{\infty}C_{t,k}>\sqrt{x})}{P(B>\sqrt{x})}=
\lim_{k_0\to\infty} \limsup_{x\to\infty}\frac{P((\sum_{k =
k_0}^{\infty}C_{t,k})^2>x)}{P(B^2>x)}=0.
$$
Repeating a similar argument as for $\alpha\in (0,1)$, we obtain
\begin{eqnarray*}
\lefteqn{\frac{P\left(\left(\sum_{k = k_0}^{\infty} C_{t,k}\right)^2>x\right)}{P(B^2>x)}}\\
& \leq & \frac{P(\sup_{k \geq k_0} B_{t-k}^2>x(1-\varepsilon)/\mu^{2k})}{P(B^2>x)}
+ \frac{P\left(\left(\sum_{k = k_0}^{\infty} C_{t,k} \1_{\{B_{t-k}^2< x(1-\varepsilon)/\mu^{2k} \}}\right)^2 >x\right)}{P(B^2>x)}.
\end{eqnarray*}
The first term could be treated using the Potter's bound since $B^2$ is regularly varying. Using Markov inequality, the second one can be
 bounded above by:
\begin{eqnarray*}
\lefteqn{
\frac{E\left(\sum_{k = k_0}^{\infty} C_{t,k} \1_{\{B_{t-k}^2< x(1-\varepsilon)/\mu^{2k} \}}\right)^2 }{xP(B^2>x)}\le
\frac{E\left(\sum_{k = k_0}^{\infty} C_{t,k}^2 \1_{\{B_{t-k}^2< x(1-\varepsilon)/\mu^{2k} \}}\right) }{xP(B^2>x)}+ }\\
& &
+\frac{E\left(\sum_{k,l = k_0\atop k\not=l}^{\infty} C_{t,k}C_{t,l} \1_{\{B_{t-k}^2< x(1-\varepsilon)/\mu^{2k}\}}\1_{\{B_{t-l}^2< x(1-\varepsilon)/\mu^{2k}  \}}\right) }{xP(B^2>x)}=:J_1(k_0)+J_2(k_0).
\end{eqnarray*}

Since by Lemma \ref{lem:bounds}:
\begin{equation*}
m_2(k) = E\left(\tilde A^{(k)}\right)^2 < C (k+1) \mu^k\,,
\end{equation*}
it follows that
\begin{eqnarray*}
E\left(\sum_{j=1}^n\tilde A_j^{(k)}\right)^2
&=&\sum_{j=1}^n E\left(\tilde A_j^{(k)}\right)^2 + \sum_{j,l=1\atop j\not=l}^n
E\left(\tilde A_j^{(k)}\right)E\left(\tilde A_l^{(k)}\right)
< Cn (k+1) \mu^{k} + n^2\mu^{2k}.
\end{eqnarray*}
Therefore,
\begin{eqnarray*}
\lefteqn{J_1(k_0)=\frac{E\left(\sum_{k = k_0}^{\infty} C_{t,k}^2 \1_{\{B_{t-k}^2< x(1-\varepsilon)/\mu^{2k} \}}\right) }{xP(B^2>x)} }\\
&\le & \sum_{k=k_0}\frac{\sum_{n\leq (x(1-\varepsilon))^{1/2}/\mu^{k} }
 E\left(\sum_{j=1}^n\tilde A_j^{(k)}\right)^2  P(B_{t-k}=n)}{xP(B^2>x)}\\
&\le &  \sum_{k=k_0}^{\infty} C(k+1) \mu^k \frac{E\left(B \1_{\{B^2<
x(1-\varepsilon)/\mu^{2k} \}}\right) }{xP(B^2>x)}
+\sum_{k=k_0}^{\infty}\mu^{2k}
\frac{E\left(B^2\1_{\{B^2< x(1-\varepsilon)/\mu^{2k} \}}\right)}{xP(B^2>x)}\\
%&\le & \sum_{k=k_0}^{\infty} C(k+1) \mu^k \frac{E\left(B \1_{\{B^2< x(1-\varepsilon)/\mu^{2k} \}}\right) }{xP(B^2>x)}+J_{12}(k_0)
&=&J_{11}(k_0)+J_{12}(k_0)\,.
\end{eqnarray*}
Since $B^2$ is regularly varying with index $\alpha/2 \in (0,1)$, Karamata's theorem applies again and
we finally have
$$
\lim_{k_0\to\infty} \limsup_{x\to\infty}J_{12}(k_0)=0.
$$
If $E(B)<\infty$, then $J_{11}(k_0)$ is bounded by
$$
J_{11}(k_0)\le C \sum_{k=k_0}^{\infty}(k+1)\mu^k \frac{E(B)}{xP(B^2>x)}
$$
and hence goes to 0 as $x\to\infty$. Otherwise, if $E(B)=\infty$
(hence $\alpha=1$),
%then we recall the following remark after Th.
%3.2 of \cite{RoSe}: If $\alpha=1$ and $E(B)=\infty$, then
%$$
%E(B 1_{\{B\le x\}})=o(x^qP(B>x)),
%$$
%for any $q>\alpha=1$. We use this with $q\in (1,2)$ to obtain
then we write
\begin{equation}\label{temp}
J_{11}(k_0)=C\sum_{k=k_0}^{\infty}(k+1)\mu^k\frac{E\left(B\1_{\{B<
\sqrt{x}\sqrt{1-\varepsilon}/\mu^{k} \}}\right)}{xP(B^2>x)}.
\end{equation}
Now, we note that
$$
E(B1_{\{B\le y\}})\le \int_0^y P(B>u)du =:\tilde L(y).
$$
Theorem 1.5.9a in \cite{BGT} yields that $\tilde L(y)$ is slowly
varying. Furthermore, the Potter's theorem yields that for any
chosen constants $A>0$ and $\delta>0$,
$$
\frac{\tilde L(y)}{\tilde L(z)}\le
A\max\left\{\left(\frac{y}{z}\right)^{\delta},\left(\frac{y}{z}\right)^{-\delta}\right\},
$$
as long as $y$ and $z$ are sufficiently large. Hence, for a
sufficiently large $x$,
$$
E\left(B\1_{\{B< \sqrt{x}\sqrt{1-\varepsilon}/\mu^{k} \}}\right)\le
\tilde L(\sqrt{x}\sqrt{1-\varepsilon}/\mu^k)\le C\tilde L(\sqrt{x})
\max\{\mu^{k \delta},\mu^{-k\delta}\}.
$$
Thus,
$$
J_{11}(k_0)\le C\frac{\tilde L(\sqrt{x})}{xP(B^2>x)}
\sum_{k=k_0}^{\infty}(k+1)\mu^{k(1-\delta)}.
$$
The series is summable if we choose $\delta\in (0,1)$. Now,
$\lim_{x\to\infty}\frac{\tilde L(\sqrt{x})}{xP(B^2>x)}=0$, since
$\tilde L$ is slowly varying and $P(B>\sqrt{x})$ is regularly
varying with index $-1/2$.
%The remark after Th. 3.2 of \cite{RoSe} gives that:
%\begin{equation}\label{eq:6b}
%\limsup_{x\to\infty}\frac{E\left(B\1_{\{B< \sqrt{x}(1-\varepsilon)/\mu^{k} \}}\right)}{x^{q/2}P(B>\sqrt{x})}<\infty\;.
%\end{equation}

Likewise,
\begin{eqnarray*}
J_2(k_0) &= & \frac{\sum_{k,l = k_0\atop k\not=l}^{\infty} E\left(C_{t,k} \1_{\{B_{t-k}^2< x(1-\varepsilon)/\mu^{2k}\}}\right) E\left(C_{t,l}\1_{\{B_{t-l}^2< x(1-\varepsilon)/\mu^{2k}  \}}\right) }{xP(B^2>x)}\\
&\le & \frac{\sum_{k,l = k_0\atop k\not=l}^{\infty} E\left((k+1)\mu^k B_{t-k} \1_{\{B_{t-k}^2< x(1-\varepsilon)/\mu^{2k}\}}\right) E\left((l+1)\mu^l B_{t-l}\1_{\{B_{t-l}^2< x(1-\varepsilon)/\mu^{2k}  \}}\right) }{xP(B^2>x)}.
\end{eqnarray*}
Again, if $E(B)<\infty$, then the term is bounded by
$$
\frac{(E(B))^2}{xP(B^2>x)}   \sum_{k,l = k_0\atop k\not=l}^{\infty}
(k+1)(l+1)\mu^{k+l}
$$
and hence
$$
\lim_{k_0\to\infty} \limsup_{x\to\infty}J_2(k_0)=0.
$$
Otherwise, if $E(B)=\infty$, we apply the same argument as for
$J_{11}(k_0)$.
%then
%$$
%\frac{\left\{E\left(B\1_{\{B^2< x(1-\varepsilon)/\mu^{2k} \}}\right)\right\}^2}{xP(B^2>x)}=
%\frac{E\left(B\1_{\{B< \sqrt{x}(1-\varepsilon)/\mu^{k} \}}\right)}{x^{q/2}P(B>\sqrt{x})}\frac{E\left(B\1_{\{B< \sqrt{x}(1-\varepsilon)/\mu^{k} \}}\right)}{x^{1-q/2}}.
%$$
%Now, $E\left(B\1_{\{B< \sqrt{x}(1-\varepsilon)/\mu^{k}
%\}}\right)=\ell(\sqrt{x})$, where $\ell$ is a slowly varying
%function. Thus, $\limsup_{x\to\infty}J_2(k_0)=0$. \hfill $\Box$

\begin{rem}\label{rem:alpha3}{\rm
Let us briefly discuss the case of $2\leq \alpha<3$ (the same
applies to all $\alpha\ge 2$). To prove Theorem \ref{eq:Theorem-1}
for this range of $\alpha$ it suffices to show that:
$$\lim_{k_0\to\infty} \limsup_{x\to\infty}\frac{P(\sum_{k =
k_0}^{\infty}C_{t,k}>x^{1/3})}{P(B> x^{1/3})}=
\lim_{k_0\to\infty} \limsup_{x\to\infty}\frac{P((\sum_{k =
k_0}^{\infty}C_{t,k})^3>x)}{P(B^3>x)}=0.
$$
Note that
\begin{eqnarray*}
\lefteqn{\frac{P\left(\left(\sum_{k = k_0}^{\infty} C_{t,k}\right)^3>x\right)}{P(B^3>x)}}\\
& \leq & \frac{P(\sup_{k \geq k_0} B_{t-k}^3>x(1-\varepsilon)/\mu^{3k})}{P(B^3>x)}
+ \frac{P\left(\left(\sum_{k = k_0}^{\infty} C_{t,k} \1_{\{B_{t-k}^3< x(1-\varepsilon)/\mu^{3k} \}}\right)^3 >x\right)}{P(B^3>x)}.
\end{eqnarray*}
The first term in this bound is handled as before.
For the second term apply Markov inequality and note that:
\begin{eqnarray*}
\lefteqn{{E\left(\sum_{k = k_0}^{\infty} C_{t,k} \1_{\{B_{t-k}^3< x(1-\varepsilon)/\mu^{3k} \}}\right)^3 } }\\
&\le &
\Big[{ \sum_{k = k_0}^{\infty} E\left( C_{t,k}^3 \1_{\{B_{t-k}^3< x(1-\varepsilon)/\mu^{3k} \}} \right)}\\
&+&
3 \sum_{k,j = k_0 \atop k\not=j }^{\infty} E\left( C_{t,k}^2 \1_{\{B_{t-k}^3< x(1-\varepsilon)/\mu^{3k} \}}
C_{t,j} \1_{\{B_{t-j}^3< x(1-\varepsilon)/\mu^{3j}\}} \right) \\
&+&
\sum_{k,j,l = k_0 \atop k\not=j\not=l }^{\infty}  E\left({
C_{t,k} \1_{\{B_{t-k}^3< x(1-\varepsilon)/\mu^{3k} \}}
C_{t,j} \1_{\{B_{t-j}^3< x(1-\varepsilon)/\mu^{3j} \}}
C_{t,l} \1_{\{B_{t-l}^3< x(1-\varepsilon)/\mu^{3l} \}} } \right)\Big] \\
& =:& L_1(k_0)+L_2(k_0)+L_3(k_0)\,.
\end{eqnarray*}
Now, $L_3(k_0)$ can be bounded above by:
\begin{equation*}
(E(B))^3  \sum_{k,j,l = k_0 \atop k\not=j\not=l
}^{\infty}(k+1)(j+1)(l+1) \mu^{k+j+l} \leq (E(B))^3 \left(
\sum_{k=k_0}^{\infty} \mu^{k} \right)^3 \,.
\end{equation*}
Hence $\lim_{k_0\to\infty} \limsup_{x\to\infty}  L_3(k_0) (xP(B^3>x))^{-1} =0 $.
To deal with $L_1(k_0)$ and $L_2(k_0)$ we proceed in the same way like in the step 3 of the proof of Theorem \ref{eq:Theorem-1}.
We only need a bound similar to that of Lemma
\ref{lem:bounds}. Using similar arguments like in the proof of Lemma \ref{lem:bounds}
observe that:
\begin{align*}
m_3(k) &\le  \mu E(\tilde A^{(k-1)})^3+3  \mu E(A^2)E(\tilde A^{(k-1)})^2+ E(A^3)\mu^{3k} \\
& =\mu \left\{\mu E(\tilde A^{(k-2)})^3+3\mu E(A^2) E(\tilde
A^{(k-2)})^2+E(A^3)\mu^{3(k-1)}\right\}\\
&\qquad +3  \mu E(A^2)E(\tilde
A^{(k-1)})^2+ E(A^3)\mu^{3k}\\
& =\mu^2 \left\{\mu E(\tilde A^{(k-3)})^3+3\mu E(A^2) E(\tilde
A^{(k-3)})^2+E(A^3)\mu^{3(k-2)}\right\}\\
&\qquad + 3\mu^2 E(A^2)E(\tilde A^{(k-2)})^2+   3  \mu
E(A^2)E(\tilde
A^{(k-1)})^2+E(A^3)\mu^{3(k-1)}+ E(A^3)\mu^{3k}\\
&= \cdots \\
&=\mu^k E(\tilde A^{(0)})^3 +3 \sum_{j=1}^{k-1} \mu^j E(A^2)E(\tilde
A^{(k-j)})^2+E(A^3)\sum_{j=1}^{k-1}\mu^{3(k-j)} .
\end{align*}
Using Lemma \ref{lem:bounds} produces:
\begin{align*}
m_3(k) &\le \mu^k E(A^3)+ E(A^3)
k\mu^{2k}+3E(A^2)E(A^2)\sum_{j=1}^{k-1}\mu^j(k-j+1)\mu^{k-j}\\
&\le \mu^k E(A^3)+ E(A^3) k\mu^{2k}+3E(A^2)E(A^2)k^2\mu^k .
\end{align*}
Similar computation can be done, in principle, for arbitrary
$m_r(k)$, $r\ge 4$. We note in passing that classical inequalities,
like Rosenthal's inequality, do not seem to be applicable here.
}
\end{rem}

\subsection{Regularly varying offspring  (Model II)}

Throughout this subsection we will assume the following conditions:
\begin{equation}\label{eq:1b}
0< \mu:=E(A)<1 \; ,
\end{equation}
\begin{equation}\label{eq:2b}
P(A>x)=x^{-\alpha}L(x) \; ,
\end{equation}
for some $\alpha\in (1,2)$ and a slowly varying function $L$.
We consider here the case $\alpha\in (1,2)$ only. For $\alpha>2$
the proof of the main result of this subsection could be adopted along
the lines of Remark \ref{rem:alpha3}.
We will also assume that the tail of $B$ is not heavier than that of $A$ in the sense that
\begin{equation}\label{eq:3b}
\lim_{x\to\infty}\frac{P(B>x)}{P(A>x)}=c\; ,
\end{equation}
where $c$ is finite (possible equal $0$) constant. If $c>0$ we need also to assume that $B$ is consistently varying.
In particular, we note that (\ref{eq:3b}) together with (\ref{eq:1b}) implies that $E(B)<\infty$.

From \cite{dima}, using induction, we can conclude that:
\begin{eqnarray*}\label{dimaas1}
\lefteqn{P(\tilde{A}^{(k)}>x)\sim \mu P(\tilde{A}^{(k-1)}>x)+ P\left(A>\frac{x}{E\tilde{A}^{(k-1)}}\right)}
\\&&\sim \mu P(\tilde{A}^{(k-1)}>x)+\mu^{(k-1)\alpha}P(A>x)\\
&& \sim \mu^2 P(\tilde{A}^{(k-2)}>x) +\mu \mu^{(k-2)\alpha}P(A>x)+\mu^{(k-1)\alpha}P(A>x)\\
&& \sim  \mu^{(k-1)\alpha}\sum_{j=0}^{k-1}\mu^{j(1-\alpha)} P(A>x)=:d_kP(A>x) .
\end{eqnarray*}
Hence, using again \cite{dima}, we get:
\begin{align}
P(C_{t,k}>x)&=P\left(\bigotimes_{i=0}^{k-1} \theta_{t-i}^{(t-k)}
(B_{t-k})>x \right) \sim E(B) P(\tilde{A}^{(k)}>x) +
P(B>x/E\tilde{A}^{(k)}) \nonumber
\\
& \sim E(B) P(\tilde{A}^{(k)}>x) + c(E\tilde{A}^{(k)})^{\alpha} P(A>x)\nonumber \\
& \sim \left(E(B) d_k+cE(\tilde{A}^{(k)})^\alpha\right)P(A>x)\nonumber\\
&=\left(E(B)d_k+c\mu^{k\alpha}\right)P(A>x)
=:\psi_kP(A>x).\label{dimaas}
\end{align}

\begin{thm}\label{eq:Theorem-2}
Under the conditions (\ref{eq:1b})-(\ref{eq:3b}), we have
\begin{equation}\label{eq:tail-behaviour2}
\lim_{x\to\infty} \frac{P(X_t>x)}{P(A>x)}
= \lim_{x\to\infty} \frac{P(X'_t>x)}{P(A>x)}
 =\sum_{k=0}^{\infty}\psi_k \; .
%\left(EB+c\mu^{k\alpha}\right)k\mu^{(k-1)\alpha}\; .
\end{equation}
\end{thm}
\noindent {\it Proof.}
We will follow the arguments of the proof of Theorem \ref{eq:Theorem-1}.
In first step one can observe from (\ref{dimaas}) that
\begin{equation*}
\lim_{x\to\infty} \frac{P(X_{t,m}>x)}{P(A>x)} =\lim_{x\to\infty}
\frac{P(X'_{t,m}>x)}{P(A>x)} =\sum_{k=0}^{m}\psi_k \;.
%=\sum_{k=0}^{m}\left(EB+c\mu^{k\alpha}\right)k\mu^{(k-1)\alpha}\;.
\end{equation*}
Step 2 is the same like i the proof of Theorem \ref{eq:Theorem-1}.
Step 3 concerns the proof of equality:
\begin{equation}\label{eq:upper_bound_Xt2}
\lim_{k_0\to\infty} \limsup_{x\to\infty}\frac{P(\sum_{k =
k_0}^{\infty}C_{t,k}>x)}{P(A>x)}=0.
\end{equation}
Similarly, one will obtain an upper bound for the process $(X'_t)$.
Now,
\begin{eqnarray*}
\lefteqn{\frac{P(\sum_{k = k_0}^{\infty} C_{t,k}>x)}{P(B>x)}} \\
& \leq & \frac{P(\sup_{k \geq k_0} B_{t-k}>x(1-\varepsilon)/\mu^k)}{P(A>x)}
+ \frac{P(\sum_{k = k_0}^{\infty} C_{t,k} \1_{\{B_{t-k}< x(1-\varepsilon)/\mu^k \}} >x)}{P(A>x)}.
\end{eqnarray*}
The first term on the right hand side is bounded above by:
$$
 \sum_{k = k_0}^{\infty} \frac{P(B>x(1-\varepsilon)/\mu^k)}{P(A>x)}\leq
(c+1)\sum_{k = k_0}^{\infty} \frac{P(A>x(1-\varepsilon)/\mu^k)}{P(A>x)}.
$$
Using the Potter's bounds we can conclude that this term goes to $0$
as first $x$ and then $k_0$ go to $\infty$. To prove that the second
term on the right hand side goes to $0$ as first $x$ and then $k_0$
go to $\infty$ we can proceed in the same like in the proof of
Theorem \ref{eq:Theorem-1} using a fact that
%$\lim_{x\to\infty}xP(A>x)=\infty$ for $0<\alpha<1$ and
$\lim_{x\to\infty}xP(A^2>x)=\infty$ for $1<\alpha<2$. \hfill $\Box$

\section{Asymptotics for sums and maxima}\label{Asymptotics for sums and maxima}

Throughout this section we assume that $(X_t)_{t\in\mathbb{Z}}$ is a
stationary process satisfying \eqref{eq:defINAR} with  distributions
of $A$ and $B$ satisfying assumptions \eqref{eq:1}--\eqref{eq:3a} of
Model I (with $\alpha\neq 1$) or \eqref{eq:1b}--\eqref{eq:3b} of
Model II. In either case, by the results of Section~\ref{Tail
behaviour}, the distribution of each $X_t$ is regularly varying with
some index $\alpha>0$.

\begin{rem}\label{Remark:mixing}
The sequence $(X_t)$ can be shown to satisfy the well known
drift/majorization criterion for geometric ergodicity of Markov
chains (cf. \cite{MeTw} or \cite{Jo04}). For $0<\varepsilon < \min
\{1,\alpha\}$ for instance, the function $V:\bbn_0 \to [1,\infty)$
given by $V(x) = 1+ x^\varepsilon$, satisfies drift condition (5) in
\cite{Jo04}, where for small set $C$ one can take a set of the form
$\{0,1,\ldots , M\}$ with integer $M$ large enough.
 This further means that process $(X_t)$ is  strongly mixing with a geometric rate
(see Theorem 2 in \cite{Jo04}).
\end{rem}

We argue next that the stationary time series $(X_t)$ is jointly regularly varying, i.e.
all its finite dimensional distributions are regularly varying.
Namely,
for a random variable $Y_0$ with Pareto distribution
$P(Y_0 > y) = y^{-\alpha}$ for $y \ge 1$ and a sequence
$Y_n = Y_0 \mu^n, \ n\in \bbn_0$, the following holds.
\begin{lem}\label{lem:JointRegVar}
Under assumptions of Model I or Model II  as $x
\to \infty$,
\begin{equation}\label{e:tailprocess}
  \bigl( (x^{-1}X_n)_{n \in \bbn_0} \, \big| \, X_0 > x \bigr)
  \fidi (Y_n)_{n \in \bbn_0}.
\end{equation}
\end{lem}
Here we use ``$\fidi$'' to denote convergence of the
finite-dimensional distributions. In the language of \cite{BaSe},
$(Y_t)$ represents the  {\em tail sequence} of the sequence $(X_t)$.

\noindent {\it Proof.}
This follows immediately, since by regular variation of $X_0$ and
the law of large numbers, as $x \toi$
\begin{eqnarray*}
  & \law( X_0/x \mid X_0 > x )
    \dto \law( Y_0 )\\
  & \law\left( {X_0}^{-1} \sum_{i=1}^{X_0} \tilde A_i^{(k)} \mid X_0 > x \right)
    \dto \delta_{\mu^k}.
\end{eqnarray*}
Hence, \eqref{e:tailprocess} follows by Slutsky argument noting that
for $t \geq 0$:
\begin{equation} \label{e:Xt}
 X_ t = \sum_{k=0}^{t-1} C_{t,k} + \sum_{i=1}^{X_0} \tilde A_i^{(t)}\,,
\end{equation}
where $C_{t,k}$'s are defined in \eqref{eq:stat-representation} and $\tilde A_i^{(t)}$'s in \eqref{e:tildeA},
with all the random variables on the right hand side being independent.
\halmos

Denote in the sequel by $(a_n)$ a sequence
of constants such that for any $u >0$ as $n \toi$
\begin{equation}\label{an}
nP(X_0 >a_n u )\to u ^{-\alpha}.
\end{equation}
It exists  and tends to $\infty$ by the regular variation property
of the random variables $X_t$. We observe next that the strong
mixing property of the process $(X_t)$ implies condition
$\mathcal{A}'(a_{n})$ of \cite{BKS}. The condition states that for
some sequence of integers $r_n\toi$, $r_n={\rm o}(n)$, denoting
$k_{n} = \lfloor n / r_{n} \rfloor$, as $n \to \infty$ we have
\begin{equation}\label{e:mixcon}
 E \biggl[ \exp \biggl\{ - \sum_{i=1}^{n} f \biggl(\frac{i}{n}, \frac{X_{i}}{a_{n}}
 \biggr) \biggr\} \biggr]
 - \prod_{k=1}^{k_{n}} \Exp \biggl[ \exp \biggl\{ - \sum_{i=1}^{r_{n}} f \biggl(\frac{kr_{n}}{n}, \frac{X_{i}}{a_{n}} \biggr) \biggr\} \biggr] \to 0.
\end{equation}
 for every $f \in C_{K}^{+}([0,1] \times
\overline{\bbr}\setminus\{0\})$, where $C^+_K(E)$ is the space of
all nonnegative continuous functions on $E$ with compact support.
Moreover, by the geometric decay of mixing coefficients (see
Remark~\ref{Remark:mixing}) such a sequence $(r_n)$ can be taken to
be of the order $o(n^\varepsilon)$, for any $\varepsilon\in (0,1)$.

Finally, the structure of regularly varying process $(X_t)$ in
either of two regimes (Model I and II)
 allow us to apply Karamata's theorem  and show that the process satisfies another well known condition in the
 literature. It is called anticlustering condition by \cite{DaHs95} and \cite{BaJaMiWi09},
 or finite mean cluster size condition in ~\cite{BKS}.
\begin{lem}\label{lem:anticl}
Under assumptions of Model I or Model II with $\alpha\not=1$ we have
\begin{equation}\label{e:anticluster}
  \lim_{m \to \infty} \limsup_{n \to \infty}
  P \biggl( \max_{m \le |t| \le r_{n}} X_{t} > ua_{n}\,\bigg|\,X_{0}>ua_{n} \biggr) = 0.
\end{equation}
\end{lem}
\noindent {\it Proof.}
By stationarity
$$
 \lim_{m \to \infty} \limsup_{n \to \infty}
  P \biggl( \max_{m \le |t| \le r_{n}} X_{t} > ua_{n}\,\bigg|\,X_{0}>ua_{n} \biggr)
  \leq
  \lim_{m \to \infty} \limsup_{n \to \infty}
   2 \sum_{t=m}^{r_n} P \biggl( X_{t} > ua_{n}\,\bigg|\,X_{0}>ua_{n} \biggr)\,.
$$
As we have seen in \eqref{e:Xt} for $t \geq 0$
$$
 X_ t = \sum_{k=0}^{t-1} C_{t,k} + \sum_{i=1}^{X_0} \tilde A_i^{(t)}\,.
$$
Since $C_{t,k}$, $t\ge 1,k\ge 0$, is independent of $X_0$, we note
that:
$$
P \left( \sum_{k=0}^{t-1} C_{t,k} > \frac{u a_n }{2} \,\bigg|\, X_0 > u a_n \right)
\leq P \left( \sum_{k=0}^{\infty} C_{t,k} > \frac{u a_n }{2} \right)= P (X_0 >ua_{n} /2 )\,.
$$
Therefore, since $r_n = o(n)$, it follows that:
\begin{equation}\label{eq:bound}
 \lim_{m \to \infty} \limsup_{n \to \infty}
 \sum_{t=m}^{r_n} P \biggl( X_{t} > ua_{n}\,\bigg|\,X_{0}>ua_{n} \biggr)
  \leq
  \lim_{m \to \infty} \limsup_{n \to \infty}
   \sum_{t=m}^{r_n} P \biggl(  \sum_{i=1}^{X_0} \tilde A_i^{(t)} > \frac{ua_{n}}{2}\,\bigg|\,X_{0}>ua_{n} \biggr)\,.
\end{equation}
However, for $\alpha > 1$, by Markov inequality:
$$
 P \biggl(  \sum_{i=1}^{X_0} \tilde A_i^{(t)} > \frac{ua_{n}}{2}\,\bigg|\,X_{0}>ua_{n} \biggr)
\leq
 2 \mu^t \dfrac{E(X_0 \1_{\{ X_0 > u a_n \}} )}{ P(X_0 > u a_n) u a_n  }.
$$
Using Karamata's theorem, note that the right hand side converges
to $2 \mu ^ t (\alpha-1)^{-1}$ as $n \toi$. Therefore,
\eqref{e:anticluster} holds for $\alpha>1$ and both Model I and II.

For Model I and $\alpha< 1$, we first observe that:
$$
P  \left(  \sum_{i=1}^{X_0} \tilde A_i^{(t)} > x \right )
\leq P \left(  {X_0}> \frac{(1-\varepsilon ) x}{\mu ^t } \right )
+
\sum_{k \leq \frac{(1-\varepsilon ) x}{\mu ^t }}  P(X_0 = k) P  \left(  \sum_{i=1}^{k} \tilde A_i^{(t)} > x \right )\,.
$$
By Markov inequality, the last sum above is bounded by
$$
\mu^t \dfrac{E(X_0 \1_{\{ X_0 \leq  \frac{(1-\varepsilon) x }{ \mu^t
} \}} )}{ x  }\,.
$$
Clearly
$$
P  \left(  \sum_{i=1}^{X_0} \tilde A_i^{(t)} >  \frac{u a_n }{2} \,\bigg|\, X_0 > u a_n  \right )
\leq
\dfrac{P  \left(  \sum_{i=1}^{X_0} \tilde A_i^{(t)} >  \frac{u a_n }{2} \right )}{P\left ( X_0 > u a_n  \right )}.
$$
Hence, the summands on the right hand side of (\ref{eq:bound}) can be bounded from
above by
$$
\dfrac{P  \left( X_0  >  \frac{(1-\varepsilon) u a_n}{2 \mu^t } \right )}{P\left ( X_0 > u a_n  \right )}
+
\mu^t    \dfrac{E(X_0 \1_{\{ X_0 \leq  \frac{(1-\varepsilon) u a_n }{2  \mu^t } \}} )}{ P(X_0 > u a_n) u a_n / 2 }\,.
$$
Now, if $\alpha<1$, apply  Potter's bound together with Karamata
theorem to conclude that after summing these terms for $t=m, \ldots
, r_n$, we may let first $n\toi$ and then $m\toi$ to obtain limit 0.
If $\alpha=1$, then we follow the same argument as in the proof of
Theorem \ref{eq:Theorem-1}. We have
$$
E(X_01_{\{X_0\le y\}})\le \int_0^y P(X_0>u)du =:\tilde L(y),
$$
where $\tilde L(y)$ is slowly varying. Hence, by Potter's bounds for
each $\delta>0$ there exists $C=C_{\delta}$ such that
\begin{align*}
\mu^t    \dfrac{E(X_0 \1_{\{ X_0 \leq  \frac{(1-\varepsilon) u a_n
}{2  \mu^t } \}} )}{ P(X_0 > u a_n) u a_n / 2 }&= \mu^t\dfrac{\tilde
L\left(\frac{(1-\varepsilon) u a_n }{2  \mu^t }\right)}{P(X_0 > u
a_n) u a_n / 2}\\
&\le C_{\delta}\mu^t\dfrac{\tilde L\left(ua_n\right)}{P(X_0
> u a_n) u a_n
/2}\max\{(2\mu^t)^{\delta}(1-\varepsilon)^{-\delta},(2\mu^t)^{-\delta}(1-\varepsilon)^{\delta}\}.
\end{align*}
\hfill $\Box$
\begin{rem}{\rm
We note that the above argument does not work for $\alpha=1$.
Indeed, Theorem 1.5.9a in \cite{BGT} implies that
$$\frac{\int_0^yu^{-1}L(u)du}{L(y)}\sim\frac{\int_0^yP(X_0>u)du}{yP(X_0>y)}\to
\infty.$$ }
\end{rem}

\subsection{Point process and maxima}

 Lemmas~\ref{lem:JointRegVar} and \ref{lem:anticl} allow us to describe the asymptotic behavior of the following point processes
\begin{equation}
\label{E:ppspacetime}
  N_{n} = \sum_{i=1}^{n} \delta_{(i / n,\,X_{i} / a_{n})} \qquad \text{ for all $n\in\bbn$.}
\end{equation}
We recall the basic notions of point processes theory, for a good
introduction we refer to \cite{Resnick87} or \cite{Re07}. Let $E$ be
a locally compact Hausdorff topological space and let $M_p(E)$ be a
space of Radon point measures on $E$. The space $M_p(E)$ is equipped
with vague metric $d(\cdot,\cdot)$. We say that a sequence $\mu_n\in
M_p(E)$ converges vaguely to $\mu\in M_p(E)$ if $\int_Efd\mu_n\to
\int_Efd\mu$ for all nonnegative continuous functions on $E$ with
compact support. Once the state $M_p(E)$ is endowed with a vague
topology, one can  study distributional limits of its random
elements like $N_n$.

It turns  out by Theorem~2.3 in \cite{BKS} that there exist a point
processes $N^{(u)}\,,\ u>0$ on the space
 $[0, 1] \times (u, \infty)$
with compound Poisson  structure such that as $n\toi$
\begin{equation} \label{e:ppconvergence}
   N_n {\bigg|_{[0, 1] \times (u, \infty)}\,} \dto N^{(u)}\,.
%    = \sum_i \sum_j \delta_{(T^{(u)}_i, u Z_{ij})} \bigg|_{[0, 1] \times \EE_u}\,.
\end{equation}

By the same theorem in \cite{BKS}  {\em the extremal index} of the
stationary sequence $(X_t)$ (see ~\cite{Resnick87} for instance)
equals
\begin{equation}\label{eq:theta}
 \theta = P (\sup_{i\geq 1} Y_i \leq 1) = P (Y_0 \leq 1/\mu) = 1 - \mu ^{\alpha}>0.
\end{equation}
Therefore, partial maxima of the process $(X_t)$ converge to the scaled Fr\'echet distribution.
\begin{cor} \label{C:maxima}
Under assumptions of Model I or Model II with $\alpha\not=1$, 
 as $n \toi $ it holds that
 $$
 P \left( \frac{M_n}{a_n} \leq x \right) \to \exp\left( - (1 - \mu ^{\alpha})
 x^{-\alpha}\right)\, ,
 $$
 for every $x\geq 0$,
 where $M_n=\max(X_1,\ldots,X_n)$.
 \end{cor}
However more precise statement on the behavior of large values can be made. For instance,
the following proposition describes the clustering of large values in the sequence
$(X_t)$ and it follows from \eqref{e:ppconvergence} by the similar token as corollary above (see e.g. \cite{LedRootz}).

\begin{prop}
Under assumptions of Model I or Model II with $\alpha\not=1$ there
is a compound Poisson process $N^\circ$ on $[0,1]$ such that
\begin{equation}\label{e:TimePro}
N_n^\circ = \dsum_{i=1}^n \delta_{\frac{i}{n}} \1_{\{ X _i > a_n
\}}
\dto N^\circ,\\
\qquad n \to \infty.
\end{equation}
\end{prop}

The proof of the proposition follows directly from
\eqref{e:ppconvergence} by Theorem 4.2 of ~\cite{Kallenberg83}.
Moreover, the limiting process $N^\circ$ has the following
representation
$$
N^\circ \eind \dsum_{i=1}^\infty \kappa_i
\delta_{T_i},
$$
where $\sum_{i} \delta_{T_i}$ is a
homogeneous Poisson point process on the interval $[0,1]$ with
intensity $\theta$ and $(\kappa_i)_{i\geq 1}$ is a sequence of
i.i.d. random variables with values in $\bbn$ independent of it.
Finally, random variables $\kappa_i$ have geometric distribution since
by Theorem 2.3 in \cite{BKS}
\begin{equation}\label{ClusSize}
P (\kappa_1=k)
=
 \frac{1}{\theta} \left[  P \left( | \{  j\geq 0 :
Y_j
> 1 \}| =k \right) - P \left( | \{  j\geq 0 : Y_j > 1 \}
| =k+1 \right) \right]
=\mu^{-\alpha (k-1)} (1-\mu^{-\alpha })
\nonumber
\end{equation}
for all $k\in \bbn$.

\subsection{Partial sums}

In the case $\alpha>2$, which we do not consider here in detail (see
Remark~\ref{rem:alpha3}), one can show that the classical central
limit theorem for strongly mixing sequences due to Ibragimov applies
(see \cite{Jo04}, Theorem 5). Namely, as $n\toi$
$$
% \frac{1}{\sqrt{n}}\left( \sum_{i=1}^n (X_i - EX_i) \right)
% =
 \frac{1}{\sqrt{n}}\left( S_n -   {n EB}/{ (1-\mu) } \right)
% \frac{1}{\sqrt{n}}\left( \sum_{i=1}^n X_i - n  EB/ (1-\mu) \right)
 \dto
 N(0,\sigma^2)\,,
$$
 where $\sigma^2=E(X_0)^2+\sum_{i=1}^\infty E(X_0X_i) <\infty$ and $S_n$ denote partial sums of the process, i.e. $S_n=
X_1+\cdots +X_n$, $n\geq 1$.
In the case $0<\alpha<2$, $X_t$'s have infinite variance and  an alternative limit theorem holds.
That is, there exists  an $\alpha$--stable random variable $\mathcal{S}_\alpha$
 such that properly centered and normalized partial sums converge to  $\mathcal{S}_\alpha$.
% $$
% \frac{1}{a_n}\left( S_n - m_n \right)\dto
% \mathcal{S}_\alpha\,,
%$$
%as $n\toi$.
%Note that for $0<\alpha <1$ one can set $m_n=0$.

Under the assumptions of Model I, for $\alpha\in(0,1)$, using
Theorem 3.1 in ~\cite{DaHs95} together with
lemmas~\ref{lem:JointRegVar} and \ref{lem:anticl} and the strong
mixing property of the stationary sequence $(X_t)$ (cf.
Remark~\ref{Remark:mixing}), one can deduce the following result
$$
\frac{S_n}{a_n} \dto \mathcal{S}_\alpha\,.
$$
Similarly,  when $\alpha\in (1,2)$   under assumptions of either
Model I or II, strictly stationary sequence $(X_t)$ satisfies
$$
\frac{S_n -  E \frac{X_{i}}{a_{n}}
  \1_{ \{ \frac{|X_{i}|}{a_{n}} \leq 1 \} }
}{a_n} \dto \mathcal{S}_\alpha\,,
$$
however, an additional technical condition is needed, see Theorem 3.1 in~\cite{DaHs95}
and condition (3.2) therein.
%that for all $\delta>0$
%\begin{equation}\label{ant}
%\lim_{u\downarrow 0}
%\limsup_{n\toi}
%P \left[\max_{k\leq n}
%    \left| \sum_{i=1}^{k} \left( \frac{X_{i}}{a_{n}}
%      \1_{ \{ \frac{|X_{i}|}{a_{n}} \leq u \} }  -  E \frac{X_{i}}{a_{n}}
%  \1_{ \{ \frac{|X_{i}|}{a_{n}} \leq u \} } \right)  \right| > \delta
%  \right] =0\,.
%\end{equation}
In either case the limit $\mathcal{S}_\alpha$ is an $\alpha$--stable
random variable with a L\'{e}vy triple $(0,\phi, r)$  identifiable
as in ~\cite{BaJaMiWi09}. Using results in ~\cite{BKS} (cf.
\cite{LR})
 it seems that with an additional effort, one can also show a functional limit theorem for
 the partial sums under similar conditions, but we do not pursue that here.

\bibliographystyle{amsalpha}

\end{document}